\documentclass[11pt]{article}

\usepackage{amssymb,latexsym,amsmath}

\begin{document}

\newcommand{\bfi}{\bfseries\itshape}

\makeatletter

\@addtoreset{figure}{section}

\def\thefigure{\thesection.\@arabic\c@figure}

\def\fps@figure{h, t}

\@addtoreset{table}{bsection}

\def\thetable{\thesection.\@arabic\c@table}

\def\fps@table{h, t}

\@addtoreset{equation}{section}

\def\theequation{\thesubsection.\arabic{equation}}

\makeatother

\newtheorem{thm}{Theorem}[section]

\newtheorem{prop}[thm]{Proposition}

\newtheorem{lema}[thm]{Lemma}

\newtheorem{cor}[thm]{Corollary}

\newtheorem{defi}[thm]{Definition}

\newtheorem{rema}[thm]{Remark}

\newtheorem{conje}[thm]{Conjecture}

\newtheorem{opp}[thm]{Open Problem}

\newtheorem{exempl}{Example}[section]

\renewcommand{\contentsname}{ }
\newenvironment{rk}{\begin{rema}  \em}{\end{rema}}
\newenvironment{exemplu}{\begin{exempl}  \em}{\end{exempl}}
\newenvironment{proof}{\begin{rema}  \em}{\end{rema}}

\newcommand{\comment}[1]{\par\noindent{\raggedright\texttt{#1}

\par\marginpar{\textsc{Comment}}}}

\newcommand{\todo}[1]{\vspace{5 mm}\par \noindent \marginpar{\textsc{ToDo}}\framebox{\begin{minipage}[c]{0.95 \textwidth}

\tt #1 \end{minipage}}\vspace{5 mm}\par}

\newcommand{\ea}{\mbox{{\bf a}}}
\newcommand{\eu}{\mbox{{\bf u}}}
\newcommand{\ueu}{\underline{\eu}}
\newcommand{\ueo}{\overline{u}}
\newcommand{\oeu}{\overline{\eu}}
\newcommand{\ew}{\mbox{{\bf w}}}
\newcommand{\ef}{\mbox{{\bf f}}}
\newcommand{\eF}{\mbox{{\bf F}}}
\newcommand{\eC}{\mbox{{\bf C}}}
\newcommand{\en}{\mbox{{\bf n}}}
\newcommand{\eT}{\mbox{{\bf T}}}
\newcommand{\eL}{\mbox{{\bf L}}}
\newcommand{\eR}{\mbox{{\bf R}}}
\newcommand{\eV}{\mbox{{\bf V}}}
\newcommand{\eU}{\mbox{{\bf U}}}
\newcommand{\ev}{\mbox{{\bf v}}}
\newcommand{\eve}{\mbox{{\bf e}}}
\newcommand{\uev}{\underline{\ev}}
\newcommand{\eY}{\mbox{{\bf Y}}}
\newcommand{\eK}{\mbox{{\bf K}}}
\newcommand{\eP}{\mbox{{\bf P}}}
\newcommand{\eS}{\mbox{{\bf S}}}
\newcommand{\eJ}{\mbox{{\bf J}}}
\newcommand{\eB}{\mbox{{\bf B}}}
\newcommand{\eH}{\mbox{{\bf H}}}
\newcommand{\leb}{\mathcal{ L}^{n}}
\newcommand{\eI}{\mathcal{ I}}
\newcommand{\eE}{\mathcal{ E}}
\newcommand{\hen}{\mathcal{H}^{n-1}}
\newcommand{\eBV}{\mbox{{\bf BV}}}
\newcommand{\eA}{\mbox{{\bf A}}}
\newcommand{\eSBV}{\mbox{{\bf SBV}}}
\newcommand{\eBD}{\mbox{{\bf BD}}}
\newcommand{\eSBD}{\mbox{{\bf SBD}}}
\newcommand{\ecs}{\mbox{{\bf X}}}
\newcommand{\eg}{\mbox{{\bf g}}}
\newcommand{\paromega}{\partial \Omega}
\newcommand{\gau}{\Gamma_{u}}
\newcommand{\gaf}{\Gamma_{f}}
\newcommand{\sig}{{\bf \sigma}}
\newcommand{\gac}{\Gamma_{\mbox{{\bf c}}}}
\newcommand{\deu}{\dot{\eu}}
\newcommand{\dueu}{\underline{\deu}}
\newcommand{\dev}{\dot{\ev}}
\newcommand{\duev}{\underline{\dev}}
\newcommand{\weak}{\rightharpoonup}
\newcommand{\weakdown}{\rightharpoondown}
\renewcommand{\contentsname}{ }

\renewcommand{\contentsname}{ }

\title{The variational complex of a diffeomorphisms group}

\author{Marius Buliga \\
 \\
IMB\\
B\^{a}timent MA \\
\'Ecole Polytechnique F\'ed\'erale de Lausanne\\
CH 1015 Lausanne, Switzerland\\
{\footnotesize Marius.Buliga@epfl.ch} \\
 \\
\and
and \\
 \\
Institute of Mathematics, Romanian Academy \\
P.O. BOX 1-764, RO 70700\\
Bucure\c sti, Romania\\
{\footnotesize Marius.Buliga@imar.ro}}

\date{This version: 13.01.2000}

\maketitle

\begin{abstract}
In this paper we propose a variational complex associated to a 
diffeomorphisms group with first order jet in a Lie group $M$. We study the 
structure of null lagrangians and we prove some fundamental properties of them, 
as well as their connection to differential invariants of the group action.
\end{abstract}


\maketitle

\section{Introduction}

A variational complex models the Lagrange formalism associated to a given 
structure. To my knowledge, the Lagrangian formalism has been developed on 
algebras, which are linear structures. In this paper I construct a variational 
complex associated to a group of diffeomorphisms. The central notion in this 
approach is the "null lagrangian" which is basically an integral invariant 
of the left action of the group on itself. 

The classical notion of "null lagrangian" means simply a lagrangian $W = W(x, \eu,
\nabla \eu)$ with the 
Euler-Lagrange equation
$$\frac{d}{dy_{j}} w(\eu, \nabla \eu ) - 
\frac{d}{dx_{i}} \frac{d}{d F_{ij}} w (\eu, \nabla \eu) = 0$$ 
satisfied for any $\eu$. This is equivalent to the identity ($\Omega \subset R^{n}$ is 
open, bounded and with smooth boundary): 
$$\int_{\Omega} W(x, \eu + \phi , \nabla (\eu + \phi) ) \mbox{ d}x  = 
\int_{\Omega} W(x, \eu  , \nabla \eu ) \mbox{ d}x$$
for any $\eu$ and for any $\phi \in C_{0}^{\infty}(\Omega)$. 

Any homogeneous null lagrangian $W = W(\nabla \eu)$ can be written as a linear 
combination of subdeterminants of $\nabla \eu$, cf. Ericksen \cite{eri} 
Ball, Currie \& Olver \cite{21} or Olver \cite{olv}. This particular structure 
of (classical) homogeneous null lagrangians leads to the formalization of 
the calculus of variations in the language of jets. Amongst the contributing papers 
we cite Tulczyjev \cite{tul}, Anderson \& Duchamp \cite{andu}, Olver \& Sivaloganathan
\cite{olsi}. More recently, the notion of variational bi-complex has been extended 
to arbitrary graded algebras (see, for example Verbovetsky \cite{ver} and the
references therein) and the theory has found applications in several domains,
connected to  the existence of variational principles associated to various problems. 

In this paper a (homogeneous) $M$ null lagrangian, where $M$ is a subgroup of 
$GL_{n}(R)$, is a function $W$ with the property: 
$$\int_{\Omega} W(\nabla (\eu . \phi)) \mbox{ d}x = \ \int_{\Omega} W(\nabla \eu) 
\mbox{ d}x$$
for any $\eu, \phi$, diffeomorphisms of $\Omega$ with compact support in $\Omega$,
such that for any $x \in \Omega$ we have $\nabla \eu (x), \nabla \phi (x) \in M$. 
By the way, the dot "." means function composition. 

It is questionable which is 
the structure of $M$ homogeneous null lagrangians and a construction of a variational 
complex should be done without previous knowledge of the form of null lagrangians. 

In this paper we propose a variational complex associated to a diffeomorphisms group 
with first order jet in a Lie group $M$. We study the structure of null lagrangians 
and we prove some fundamental properties of them.

\section{Preliminaries}

\subsection{Notations}

$GL_{n}(R) \subset {\mathbb R}^{n \times n}$ is the
multiplicative group of  all 
invertible, orientation preserving,
 matrices, i.e the set of all $\eF$ such that  
$det \ \eF \ > 0$.

\subsection{Basic definitions and properties}

\begin{defi}
For any Lie subgroup $M \leq GL_{n}$ we define the associated local group
\begin{equation}
[M] = \left\{ \phi \in C^{\infty}(R^{n},R^{n}) \mid  \forall \ x \in R^{n} \ \ 
\nabla \phi(x) \in M \right\}
\label{dem1}
\end{equation}
and it's subgroup of compactly supported diffeomorphisms
\begin{equation}
[M]_{c} = \left\{ \phi \in C^{\infty}(R^{n},R^{n}) \mid  \forall \ x \in R^{n} 
\nabla \phi(x) \in M \ , \ supp \ (\phi - \ id) \subset \subset R^{n} \right\}
\label{dem}
\end{equation}
\end{defi}

$\mathcal{ A}$ is the group of affine homothety-translations. Any 
element of $\mathcal{ A}$ has the form:  
$$ \alpha(x_{0},y_{0},\epsilon)(x) \ = \ \ef (x) \ = \ x_{1} + 
\epsilon (x -x_{0}) \ \ , \ \ x_{0}, 
x_{1}  \  \in \ {\mathbb R}^{n} \ , \ \epsilon \ > \ 0 \ \ .$$ 
We consider on $\mathcal{ A}$ the punctual convergence of functions 
defined on  ${\mathbb R}^{n}$ with values in ${\mathbb R}^{n}$.

 $[M]_{c}$ is a set  of functions from  
${\mathbb R}^{n}$ to  ${\mathbb R}^{n}$, 
which satisfies the following axioms:
\begin{enumerate}
\item[A1/] $([M]_{c},.)$ is a group with  
the function composition operation ".";
\item[A2/] the following action is well defined:
$$A: \mathcal{ A} \times [M]_{c} \rightarrow [M]_{c} \ \ , \ \ A( \ef, \phi )
 \ = \ \ \ef . \phi . \ef^{-1} \ \ . $$ 
\end{enumerate}

\begin{rk}
These axioms are respected also by $[M]$.
\end{rk}

\begin{defi} For any open set $E \subset {\mathbb R}^{n}$ 
we define 
$$[M](E) \ = \ \left\{ \phi \in [M] \mbox{ : } supp \ (\phi - id) 
\subset \subset E \right\} \ \ . $$

For any $x_{0} \in {\mathbb R}^{n}$  
the first order jet of $M$ in $x_{0}$ is: 
$$J^{1}(x_{0},M) \ = \ \left\{  
\nabla \phi (x_{0}) \mbox{ : } \phi \in [M] \right\} \ \ . $$
\label{dunu}
The first order jet of $M$ compactly supported diffeomorphisms is: 
$$J^{1}_{c}(x_{0},M) \ = \ \left\{  
\nabla \phi (x_{0}) \mbox{ : } \phi \in [M]_{c} \right\} \ \ . $$
\end{defi}

\begin{prop}
If $[M]_{c}$ acts transitively on $R^{n}$ then 
there is a sub semigroup $J_{c}(M)$  of the multiplicative
 group $GL_{n}(R)$ 
such that for any 
$x_{0} \in {\mathbb R}^{n}$ we have $J^{1}_{c}(x_{0},M) \ = \ 
J_{c}(M)$.  We have also $J(x_{0}, M) = M$. 
\label{p3}
\end{prop}

\paragraph{Proof.} We first prove that $J^{1}_{c}(x_{0},M)$ is 
semigroup. We have $id \in [M]_{c}$, 
hence $I$, the identity matrix, belongs to $J^{1}_{c}(x_{0},M)$. 
Let us consider 
$R, S \in J^{1}_{c}(x_{0},M)$ and  
 $\phi, \psi \in [M]_{c}$ such that $R = 
\nabla \phi (x_{0})$, 
$S  =  \nabla \psi(x_{0})$. We define the translation 
$\ef \in \mathcal{ A}$: 
$\ef(x)  =  x + \psi(x_{0}) - x_{0}$. From A2/ we have 
$ \tilde{\phi} \ = \ \ef . \phi . \ef^{-1} \in [M]_{c}$, hence  from  
$\nabla \tilde{\phi} (\psi(x_{0}))  =  \nabla \phi (x_{0})$ 
and A2/  we infer that  
$$RS \ = \ \nabla \phi (x_{0}) \nabla \psi(x_{0}) \ = \ 
\nabla \tilde{\phi} (\psi(x_{0})) 
 \nabla \psi(x_{0}) \ = $$ $$ = \ \nabla ( \tilde{\phi} . 
\psi) (x_{0})  \ \in J^{1}_{c}(x_{0}, M) \ \ . $$
A simple argument based on A2/ shows that $J^{1}(x_{0},M)$ 
does not depend on  
 $x_{0} \in {\mathbb R}^{n}$. For a fixed, arbitrarily 
chosen,  $x_{0}$ we define $J_{c}(G) \ = \ 
J^{1}_{c}(x_{0},M)$. 

The proof of the fact that if  $\eF 
\in J_{c}(G)$ then  
$\eF^{-1}$ exists and $\eF^{-1} \in J_{c}(G)$ is similar.

The last equality is trivial because we have $J(M) \subset M$ and for 
any $\eF \in M$ we have $\eF \in [M]$ (as a linear function), therefore 
$\eF \in J(G)$.  
\quad $\blacksquare$

\begin{rk}
We notice that the groups $[M]$  are 
determined by $J(M)$, that is~:
if  $J(M_{1})  = 
J(M_{2})$ then $[M_{1}]  =  [M_{2}]$.  This property justifies 
the name "local group" for a diffeomorphisms group $[M]$. 
\end{rk}

\begin{prop}
$J_{c}(M)$ is a normal Lie subgroup of $M$.
\end{prop}

\paragraph{Proof.}
The fact that $J_{c}(M)$ is a Lie group follows from the continuity properties 
of the composition between smooth diffeomorphisms. 

Consider now $\eF \in M$. Then, for any $\phi \in [M]_{c}$ we have 
$\eF . \phi . \eF^{-1} \in [M]_{c}$. We deduce that $J_{c}(M)$ is a normal 
subgroup of $M$.
\quad $\blacksquare$

\begin{exemplu}
We obviously have $[GL_{n}]_{c} \ = \ Diff^{\infty}_{0}$, 
and (not so obvious)  $J_{c}(GL_{n}) =  GL_{n}$.   
 Also, for any $\Omega \subset R^{n}$ with 
lipschitz boundary the group $[GL_{n}](\Omega)$ acts transitively on $\Omega$.
\label{ex1}
\end{exemplu}
\begin{exemplu} 
Let us consider $Diff^{\infty}_{0}(dx)$, the 
subgroup of $Diff^{\infty}_{0}$ containing all 
volume preserving smooth diffeomorphisms with compact support. 
We have $Diff^{\infty}_{0}(dx)  =  [SL_{n}(R)]_{c}$ and 
$J_{c}(SL_{n}) = SL_{n}$. For any $\Omega \subset R^{n}$ with 
lipschitz boundary the group $[SL_{n}](\Omega)$ acts transitively on $\Omega$. 
\label{ex2}
\end{exemplu}

\begin{exemplu} 
For any  $\eu: {\mathbb R}^{2n} \rightarrow {\mathbb R}^{2n}$ 
and $\omega$,the canonical symplectic 2-form on 
${\mathbb R}^{2n}$,  we 
denote by 
 $\eu^{*}(\omega)$  the transport of $\omega$. Let us define 
$Diff^{\infty}_{0}(\omega)$: 
$$Diff^{\infty}_{0}(\omega) \ = \ \left\{ \phi \in 
Diff^{\infty}_{0} \mbox{ : } 
\phi^{*}(\omega) \ = \ \omega \ \right\} \ \ .$$
 We  have the equalities:  
$$J(Diff^{\infty}_{0}(\omega) ) \ = \ Sp_{n}(R) \ = \  
\left\{ \eF \in {\mathbb R}^{2n \times 2n} \mbox{ : }  
\eF \omega \eF^{T} \ =  \ \omega \ \right\} $$
and $J_{c}(Sp_{n}) = Sp_{n}$. Also, as in the previous examples, 
 $[Sp_{n}](\Omega)$ acts transitively on $\Omega$. 
\label{ex3}
\end{exemplu}

\begin{exemplu}
Let us take $C_{n} = \left\{ \lambda \eR \mid \lambda > 0 \ ,  \ \eR \eR^{T} = I_{n} 
\right\}$, the group of linear conformal matrices. Then 
$[C_{n}]$ is the group of conformal diffeomorphisms and $[C_{n}]_{c} = \left\{ id 
\right\}$, therefore $J_{c}(C_{n}) = \left\{ I_{n} \right\}$. 
\label{ex4}
\end{exemplu}

For all the transitivity results needed in these examples we 
refer to Michor \& Vizman \cite{8}.

In the following lemmas we collect some elementary 
facts connected  to the algebraic 
structure previously introduced.

\begin{lema}
Let $A, B$ be non empty open subsets of ${\mathbb R}^{n}$.
 If $A$ is bounded then there exists  
$\ef \in \mathcal{ A}$ such that the application $A(\ef, \cdot) :
 [M](A) \rightarrow [M](B)$ is well defined, 
injective and continuous. 
\label{p1}
\end{lema}

\begin{lema}
If $A$ and $B$ are two open disjoint sets then for any  
$\phi \in [M](A)$, $\psi \in 
[M](B)$ we have $\phi . \psi \ = \ \psi . \phi \ \in \ 
[M](A \cup B)$ . 
\label{p2}
\end{lema}

\begin{lema}
For any group $[M]_{c}$ we have  $\eF.\phi. \eF^{-1} \in [M]_{c}$ if $\eF \in M$. 
\label{dcompl}
\end{lema}

\section{Null lagrangians and group invariants}

We shall denote by $T[M](\Omega)$ the space of all vector fields $\eta$ over $\Omega$ 
with the associated one-parameter flow in $[M](\Omega)$.

We associate to  any $W \in \Lambda^{s}(M)$ the integral 
$I_{W}$ defined over $[M](\Omega) \times (T[M](\Omega))^{s}$ by the formula:
$$ I_{W}(\phi, \eta) = \int_{\Omega} W(\phi, \nabla \phi) (\nabla \eta^{1}, \nabla 
\eta^{2}, 
..., \nabla \eta^{s}) \mbox{ d}x $$

$\Lambda(M)$ is the graded algebra of differential forms over $J_{c}(M)$. 
The Lie algebra associated to $J_{c}(M)$ will be denoted by $j_{c}(M)$. 
The space of constant differential forms of order $s$ over $j_{c}(M)$ is 
$\lambda^{s}(M)$. Then 
$$\Lambda^{s}(M) \ = \  C^{\infty}(R^{n} \times J_{c}(M), \lambda^{s}(M))$$

\begin{defi}
A differential form $\alpha \in \Lambda^{r}(R^{n})$ is a M differential 
invariant if for any $\phi \in [M]_{c}$ the form $\phi^{*} \alpha \ - \alpha$ 
is exact. 
\label{ddinv}
\end{defi}

One can extend this definition to "G differential invariant" simply by replacing 
"$\phi \in [M]_{c}$" with " $\phi \in G$, where $G$ is a diffeomorphisms group.

\begin{defi}
A function $W \in \Lambda^{0}(M)$ is a M null lagrangian if 
for any $\ef \in \mathcal{A}$, with the notation 
$(W, \ef)(x, \eF) = W(\ef(x), \eF)$,  we have: $I_{(W, \ef)}( \cdot ; \Omega)$ 
is constant over $[M](\Omega)$. 
\label{dnullag}
\end{defi}

\begin{prop}
The definition \ref{dnullag}  does not depend on the choice of $\Omega$ in the class 
of smooth open bounded sets.  
\label{pindep}
\end{prop}

\paragraph{Proof.}
Let $\Omega'$ smooth, bounded and $\phi' \in [M](\Omega')$. There is an 
element $\ef \in \mathcal{A}$, $\ef(x) = a + \varepsilon x$, 
such that $\ef(\Omega') \subset \Omega$. Then 
we have $\ef . \phi' . \ef^{-1} \in [M](\Omega)$. $W$ verifies definition 
\ref{dnullag}, therefore: 
$$\int_{\Omega} W(x, I) \mbox{ d}x = \int_{\Omega} W(\ef (\phi'(\ef^{-1}(x))), 
\nabla \phi'(\ef^{-1}(x))) \mbox{ d}x $$
We change variables $y = \ef^{-1} (x)$ and we obtain: 
$$\int_{\ef^{-1}(\Omega)} W(\ef (x), I) \mbox{ d}y =  \int_{\ef^{-1}(\Omega)} 
W(\ef(\phi'(y)), \nabla \phi'(y)) \mbox{ d}y$$
Because $\phi'$ has compact support in $\Omega'$, we have: 
$$\int_{\Omega'} W(\ef (y), I) \mbox{ d}y =  \int_{\Omega'} 
W(\ef(\phi'(y)), \nabla \phi'(y)) \mbox{ d}y$$
This resumes the proof. 
\quad $\blacksquare$

The class of null lagrangians which generate null integrals is denoted by 
$$ NL^{0}(M) \ = \ \left\{ W \in \Lambda^{0}(M) \ : \  I_{W} = 0  \right\}$$
We have fixed the constant value of $I_{W}$  to be $0$ for further technical reasons. 
However this not causes problems, because we may think instead that a factorization by 
$R$ has been performed. 

The class of homogeneous null lagrangians is made by all $W = W(\nabla \eu)$ which are
null lagrangians. This class is denoted by $nl^{0}(M)$. No factorization by $R$ has 
been done in this case.

To any differential form $\alpha \in \Lambda^{r}(R^{n})$ and any 
$(v_{1}, ... , v_{r}) \in R^{nr}$ we associate a potential 
$\alpha^{*} (v_{1}, ... , v_{r}) \in \Lambda^{0}(GL_{n})$ in the following way: 
$$\alpha^{*}(y, \eF)(v_{1}, ... , v_{r}) = \alpha(y)(\eF v_{1}, ... , \eF v_{r}) $$
In order to shorten the notation we shall denote by $I_{\alpha}( \cdot ; \Omega)$ 
the mapping defined over $X_{1}^{r}(M)$ with values in $R^{nr}$, given by: 
$$I_{\alpha}(\eu ; \Omega)(v_{1}, ... , v_{r}) = \int_{\Omega} 
\alpha^{*}(\eu(x), \nabla \eu (x))(v_{1}, ... , v_{r}) \mbox{ d}x $$

\begin{prop}
If $\alpha$ is a M (or $[M](\Omega)$) differential invariant then for any 
$(v_{1}, ... , v_{r})$ the mapping $(y , \eF) \mapsto \alpha^{*}(y, \eF)
(v_{1}, ... , v_{r})$ is a null lagrangian (shortly:  $\alpha^{*}$ is a 
vectorial null lagrangian). 
\label{pdinv}
\end{prop}

\paragraph{Proof.}
We remark that generally, if $\alpha$ is exact and with compact support 
in $\Omega$, then $\alpha = d \beta$, with $\beta$ with compact support in $\Omega$. 
For any $\phi \in [M](\Omega)$ the differential form $\phi^{*} \alpha -
\alpha$ has compact support in  $\Omega$. 

We apply the definition of $\alpha^{*}$, definition \ref{ddinv} 
and integration by parts (Gauss formula) for the integral $I_{\alpha}(\eu ; \Omega)
(v_{1}, ... , v_{r})$ and we conclude the proof. 
\quad $\blacksquare$

The following theorem will be essential in further proofs. 

\begin{thm}
If $W: \Omega \times \Omega \times J_{c}(M) \rightarrow R$ is continuous and 
$$I_{W}(\eu) \ = \int_{\Omega} W(x, \eu(x), \nabla \eu(x)) \mbox{ d}x$$
is constant over $[M](\Omega)$ then for any $x_{1} \in \Omega$ and any  
$\psi \in [M](\Omega)$ the mapping 
$\eF \in J_{c}(M) \mapsto W(x_{1}, \psi(x_{1}), \nabla \psi(x_{1}) \eF)$
is a homogeneous null lagrangian. 
\label{tfundam}
\end{thm}

\paragraph{Proof.}
Let us consider $x_{1} \in \Omega$ and $h > 0$. 
 $Q_{h}$ is the cube $x_{1}^{i} < x^{i} < x^{i} + 1/h$. 
We take $\phi \in [M](Q_{1})$ and $k \in N$. The extension of 
$\phi$ by periodicity over ${\mathbb R}^{n}$ 
is denoted by $\tilde{\phi}$. 
We define then: 
$$\phi_{h,k} (x) \ =  \ \left\{ \begin{array}{ll} 
( hk)^{-1} \left( \tilde{\phi}(hk(x-x_{1}) + x_{1})  -x_{1} 
\right) 
+ x_{1} & \mbox{ if } x \in 
Q_{h} \\ 
x & \mbox{ otherwise} \ \ . 
\end{array}
\right. $$
We have  
$\phi_{h,k} \in [M](\Omega)$. Any set   $Q_{h}$  decomposes in  
$k^{n}$ cubes which will 
be denoted by $Q_{hk,j}$, $j = 1, ..., k^{n}$, such that 
$Q_{hk,1} = Q_{hk}$.  The corner of 
 $Q_{hk,j}$ with least distance from $x_{1}$ is denoted by 
$x_{j}$.  

Let us now consider  $\psi \in [M](\Omega)$. 
For a sufficiently large  
$h$ we have $Q_{h} \subset \Omega$, hence 
$I(\psi . \phi_{h,k} ; \Omega)$ makes sense. 
We decompose this integral in two parts: 
\begin{equation}
I(\psi . \phi_{h,k} ; \Omega) \ = \ I(\psi . \phi_{h,k} ; 
Q_{h}) \ + \ I(\psi  ; \Omega \setminus Q_{h}) \ \ , 
\label{prt21}
\end{equation}
$$I(\psi . \phi_{h,k} ; Q_{h}) \ =   
\sum^{k^{n}}_{j=1} \int_{Q_{hk,j}} \left[ 
W(x, \psi . \phi_{h,k}(x), \nabla ( \psi . \phi_{h,k}) (x)) 
\right. $$
\begin{equation}
  \left.  - \ W(x_{j},\psi . \phi_{h,k}(x_{j}), 
\nabla  \psi ( \phi_{h,k} (x_{j})) \nabla \phi_{h,k} (x)) 
\right] \mbox{ d}x \ + 
\label{prt22}
\end{equation}
$$  + \  \sum^{k^{n}}_{j=1} \int_{Q_{hk,j}} 
W(x_{j},\psi . \phi_{h,k}(x_{j}), 
\nabla  \psi ( \phi_{h,k} (x_{j})) \nabla \phi_{h,k} (x)) 
\mbox{ d}x \ \ . $$
Notice that $\phi_{h,k}$ converges weakly to  $id$. Because  
$W$  and $\nabla \psi$ are continuous and $\phi_{h,k}$
 converges uniformly to  $id$ with $k \rightarrow \infty$, it follows that 
the first sum from the right-handed member of the 
equality  (\ref{prt22}) converges to zero.

By the change of variable $y = hk(x-x_{j}) + x_{1}$
 we obtain: 
\begin{equation}
 \int_{Q_{hk,j}} W(x_{j},\psi . \phi_{h,k}(x_{j}), 
\nabla  \psi ( \phi_{h,k} (x_{j})) \nabla \phi_{h,k} (x))  
\mbox{ d}x \ = 
\label{key1}
\end{equation}
$$ = \ (hk)^{-n} \ 
\int_{Q_{1}} W(x_{j}, \psi(x_{j}), \nabla \psi (x_{j})
 \nabla \phi (y) ) \mbox{ d}y \ \ . $$
We deduce from here that the second sum of the right-handed
 member  (\ref{prt22}) is a Cauchy sum. By a passage
 to the limit as $k \rightarrow \infty$  we get the equality:   
\begin{equation}
\label{prt23}
\lim_{k \rightarrow \infty} I(\psi . \phi_{h,k} ; Q_{h}) \ = \ 
\int_{Q_{h}} \int_{Q_{1}} W(x, \psi(x), \nabla \psi(x) 
\nabla \phi(y)) \mbox{ d}y \mbox{ d}x \ \ . 
\end{equation}

From  the hypothesis we have: 
$$ \lim_{k \rightarrow \infty} I(\psi . \phi_{h,k} ; \Omega) 
 = \ \lim_{k \rightarrow \infty} 
I(\psi . \phi_{h,k} ; Q_{h})      + 
\ I(\psi  ; \Omega \setminus Q_{h}) $$   
\begin{equation}
 = \ I(\psi ; Q_{h}) \ + \ I(\psi  ; \Omega \setminus Q_{h}) 
\ \ , 
\end{equation} 
therefore (\ref{prt23}) implies that: 
\begin{equation} \label{prt24} 
\int_{Q_{h}} \int_{Q_{1}} W(x, \psi(x), \nabla \psi(x) \nabla 
\phi(y)) \mbox{ d}y \mbox{ d}x  \ = 
\end{equation}
$$= \ \int_{Q_{h}} W(x, \psi(x), \nabla \psi(x) ) \mbox{ d}x 
\ \ .$$ 
We multiply the relation (\ref{prt24}) with $h^{n}$ and  pass to 
the limit as  $h \rightarrow \infty$. 
The result is: 
$$\int_{Q_{1}} W(x_{1}, \psi(x_{1}), \nabla \psi (x_{1}) \nabla 
\phi (y)) \mbox{ d}y \ 
= \ W(x_{1}, \psi(x_{1}) , \nabla \psi (x_{1})) $$ 
which concludes  the proof.
\quad $\blacksquare$

\section{The variational complex of the diffeomorphism group}

We introduce  two  graded vector spaces:
$$K^{0}(M) = NL^{0}(M) \ , \ K^{p}(M) = \left\{ W \in \Lambda^{p}(M) \ : \ \forall H 
\in j_{c}(M), W(\cdot) (H, ...) \in K^{p-1}(M) \right\}$$
$$NL^{p}(M) = \left\{ W \in K^{p}(M) \ : \ I_{W} = 0  \right\}$$

The  graded derivation operator $D$ on $\Lambda(M)$ is defined further. 
It is sufficient to define $D^{0}$, $D^{1}$. 

$D^{0}: \Lambda^{0}(M) \rightarrow \Lambda^{1}(M)$ is defined by the formula: 
$$Dw (y,\eF) H \ = \ w(y, \eF) \mbox{ tr} H \ - \ 
\langle \frac{\partial w}{\partial \eF}(y, \eF) , \eF H \rangle$$
$D^{1}: \Lambda^{1}(M) \rightarrow \Lambda^{2}(M)$ is defined by: 
$$Dw(y, \eF) (H, P) \ = \  D(w_{ij} H_{ij}) (y, \eF) P  -  D(w_{ij} P_{ij}) (y, \eF) H 
+ w(\eF) [P,H] $$
We have then the obvious proposition: 
\begin{prop}
The sequence 
$$\Lambda^{0}(M) \stackrel{D}{\longrightarrow} \Lambda^{1}(M) \stackrel{D}{\longrightarrow} 
\Lambda^{2}(M)$$ 
is semi-exact, that is $D^{2} = 0$. 
\label{pexa}
\end{prop}

\paragraph{Proof.}
We leave to he reader to verify, by direct calculus. 
\quad $\blacksquare$

Remark that if $j(M) \subset sl_{n}$ (the algebra of null trace matrices), then 
the complex defined above is simply the de Rham complex.

More interesting is the following:

\begin{prop}
The sequences 
$$NL^{p-1}(M) \stackrel{D}{\longrightarrow} NL^{p}(M) \stackrel{D}{\longrightarrow} 
NL^{p+1}(M)$$ 
$$K^{p-1}(M) \stackrel{D}{\longrightarrow} K^{p}(M) \stackrel{D}{\longrightarrow} 
K^{p+1}(M)$$ 
are well defined and, off course, semi-exact. 
\label{pnlex}
\end{prop}

\paragraph{Proof.}
We shall prove only that 
$$K^{0}(M) \stackrel{D}{\longrightarrow} K^{1}(M)$$ 
is well defined. The rest follows.  

Let $\eu \in [M](\Omega)$ and $\eta \in T[M](\Omega)$. The flow on $\Omega$ 
generated by $\eta$ is $t \mapsto \phi_{t} \in [M](\Omega)$. 

We consider the function:
$$I(\eu, t) \ = \ \int_{\Omega} w( \eu . \phi_{t}^{-1}, 
\nabla ( \eu . \phi_{t}^{-1}) ) \mbox{ d}x$$ 
If $w \in NL^{0}(M)$ then $I(\eu, t)$ is constant with respect to $t$, therefore 
the derivative of $I(\eu, t)$ relative to $t$ at $t=0$ is null, for any $\eu$. 
We have the following representation of this derivative:
$$\frac{\partial I}{\partial t}(\eu, 0) \ = \ \int_{\Omega} 
\hat{w}(x,\eu(x), \nabla \eu(x)) \mbox{ d}x$$ 
where the potential $\hat{w}$ is given by the expression: 
$$\hat{w} (x,y, \eF) =  w(y,F) \mbox{ div} \eta (x) \ - \ 
\langle \frac{\partial w}{\partial \eF}(y,\eF) , \eF \nabla \eta (x) 
\rangle$$
We deduce that $\hat{w}$ is a non-homogeneous null lagrangian. 
From  theorem \ref{tfundam} we deduce that 
$Dw$, previously defined, indeed belongs to $K^{1}(G)$. 
\quad $\blacksquare$

\begin{defi}
The complex $NL(M)$ introduced in  proposition \ref{pnlex} is called 
the complex of [M] null lagrangians. Analogously, the complex of homogeneous 
null lagrangians is $nl(M)$. 
\label{dcomnl}
\end{defi}

In the context of this paper  we introduce the following definition of the
Euler-Lagrange operator: 

\begin{defi}
The Euler-Lagrange operator associated to $w \in \Lambda^{0}(M)$ is the function: 
$$Ew: [M](\Omega) \rightarrow  T[M](\Omega)^{*}$$ 
defined by 
\begin{equation}
Ew(\eu) \eta \ = \ \int_{\Omega} \left( \frac{d}{dy_{j}} w(\eu, \nabla \eu ) - 
\frac{d}{dx_{i}} \frac{d}{d F_{ij}} w (\eu, \nabla \eu) \right) \nabla \eu_{jk}
\eta_{k} \mbox{ d}x
\label{eula}
\end{equation}
\label{deula}
\end{defi}

From the proof of the previous proposition we extract the formula:
\begin{equation}
\frac{d}{dt} \int_{\Omega} w( \eu . \phi_{t}^{-1}, 
\nabla ( \eu . \phi_{t}^{-1}) ) \mbox{ d}x \ = \ 
\int_{\Omega} Dw(\eu, \nabla \eu) \nabla \eta \mbox{ d}x 
\label{der}
\end{equation}
We integrate by parts, taking into account that $supp \ \eta \subset \subset \Omega$, 
and we obtain: 
\begin{equation}
\frac{d}{dt} \int_{\Omega} w( \eu . \phi_{t}^{-1}, 
\nabla ( \eu . \phi_{t}^{-1}) ) \mbox{ d}x \ = \ 
\int_{\Omega} \mathcal{E}_{j} w(\eu, \nabla \eu)  \nabla \eu_{jk}
\eta_{k}\mbox{ d}x
\label{der2}
\end{equation}
where $\mathcal{E}w$ is the classical Euler-Lagrange operator defined by: 
\begin{equation}
\mathcal{E}_{j} w(\eu, \nabla \eu) \ \ = \ \frac{d}{dy_{j}} w(\eu, \nabla \eu ) - 
\frac{d}{dx_{i}} \frac{d}{d F_{ij}} w (\eu, \nabla \eu) 
\label{elo}
\end{equation}

We construct next the variational complex associated to $[M]$.
The following proposition is straightforward, due to proposition \ref{pnlex}:

\begin{prop}
The following sequence is semi-exact:
\begin{equation}
\frac{\Lambda^{p}(M)}{NL^{p}(M)} \stackrel{D}{\rightarrow}
\frac{\Lambda^{p+1}(M)}{NL^{p+1}(M)} 
 \ \ ,  \ D \left( W + NL^{p}(M) \right) = D W + NL^{p+1}(M)
\label{varc}
\end{equation}
\label{pvarc}
\end{prop}

\begin{defi}
With the notation 
$$V^{p}(M) \ = \ \frac{\Lambda^{p}(M)}{NL^{p}(M)}$$
the complex $V(M)$ is called the variational complex associated to $M$. 
For any $w \in \Lambda^{p}(M)$, we shall denote by $NL w$ the equivalence 
class: 
$$NL w = w + NL^{p}(M) \ \in \  V^{p}(M) $$
\label{dvarc}
\end{defi} 

The name "variational complex" given to $V(M)$ is justified by the following 
proposition, which tells that the Euler-Lagrange operator $E w$ can be identified 
with $D \ NL w $, that is the derivation $D: V^{0} \rightarrow V^{1}$ "is" the
Euler-Lagrange operator: 

\begin{prop}
For any $w \in \Lambda^{0}(M)$, $Ew = 0 \in T[M](\Omega)^{*}$  if and only 
if $D \ NL w \ = \ 0 \in V^{1}(M)$. 
\label{jvarc}
\end{prop}

\paragraph{Proof.}
Recall the formula \eqref{der} (using also definition \ref{deula}): 
\begin{equation}
\int_{\Omega}  Dw(\eu, \nabla \eu) \nabla \eta \mbox{ d}x = \ E w (\eu) \eta 
\label{ed}
\end{equation}
If $Ew = 0 \in T[M](\Omega)^{*}$ then $w$ is constant along any one parameter flow 
in $[M](\Omega)$, hence $w \in NL^{0}$, which is the same as 
$D \ NL w \ = \ 0 \in V^{1}(M)$. Conversely, suppose that $w \in NL^{0}$. Then we 
have immediately $E w (\eu) \eta = 0$, for any $\eu \in [M](\Omega)$ and 
$\eta \in T[M](\Omega)$. 
\quad $\blacksquare$

Let us introduce another complex, made by variational integrals: 

\begin{defi}
The graded vector space of variational integrals $VI(M)$ is: 
\begin{equation}
VI^{p}(M) \ = \ \left\{ I_{W} \ : \ W \in \Lambda^{p}(M) \right\} 
\label{vid}
\end{equation}
We consider the derivation operator $D I_{W} \ = \ I_{DW}$. 
\label{ved}
\end{defi}

We have therefore, the following 
easy consequence of proposition \ref{pvarc}: 

\begin{prop}
The sequence $$VI^{p-1}(M) \stackrel{D}{\rightarrow} VI^{p}(M) 
\stackrel{D}{\rightarrow} VI^{p+1}(M)$$ 
is semi-exact. Moreover, the mapping 
$$NL \ w \in V^{p}(M) \ \mapsto \ I_{W} \in VI^{p}(M)$$ 
is an  isomorphism from the variational complex $V(M)$ to the variational 
integrals complex $VI(M)$. 
\label{pvi}
\end{prop}

\paragraph{Proof.}
Recall that we have demanded that $w \in NL^{p}(M)$ if and only if  
$I_{w} = 0$. Therefore we have $w - w' \in NL^{p}(M)$ if and only if $I_{w} = I_{w'}$. 

 The second part of the proposition is  straightforward from definition \ref{ved} 
and proposition \ref{pnlex}. 
\quad $\blacksquare$

As a conclusion, we get a more precise image of the variational 
complex if we look at  the complex of variational integrals $VI(M)$. 
It is useful to consider also the  relation, coming from \eqref{der}, 
$$ D I_{w} (\eu, \eta) = \frac{d}{dt} I_{w} (\eu . \phi_{-t})_{|_{t=0}}$$
where $w \in \Lambda^{0}(M)$ and $\phi_{t}$ is the one-parameter flow generated 
by $\eta$.

\section{Examples}

It is visible now that a central object in the construction presented in 
this paper is the space of homogeneous null lagrangians $nl^{0}(M)$. 
We shall look to this space in the followings.

We derive first a necessary condition for $W$ to be a homogeneous $M$ null lagrangian.

\begin{thm}
Let $W \ = \ W(\eF)$ be a homogeneous $M$ null lagrangian. Then for any 
$\eF \in J_{c}(G)$ and for any $\eta \in T[M](\Omega)$ we have the inequality: 
\begin{equation}
\int_{\Omega}D \left( D W(\eF) \nabla\eta(x)\right) \nabla \eta(x)
\mbox{ d}x \ = \ 0 \ \ \ .
\label{legh}
\end{equation}
\label{tlegh}
\end{thm}

\paragraph{Proof.} 
 Consider the function 
$$I(t) \ =  \ \int_{\Omega} W(\eF \nabla \phi_{-t}(x)) 
\mbox{ d}x \ \ . $$
This is a $C^{2}$ function which is constant, 
according to hypothesis upon $W$  and theorem \ref{tfundam} .  This fact implies 
that 
$$\frac{\partial I}{\partial t}(0) \ = \ 0  \ \ \ , \ \ 
\frac{\partial^{2} I}{\partial t^{2}} (0) \ = \ 0 \ \ . $$
The first relation is trivially satisfied. 

We apply twice an  integration by parts argument to the second equality  in order to
obtain the desired inequality.   
\quad $\blacksquare$

In order to see what \eqref{legh} means let us take $M = GL_{n}$. In this case we have 
$$T[M](\Omega) \ = \ C^{\infty}_{0}(\Omega, {\mathbb R}^{n})$$ 
 hence for any $\eta \in T[M](\Omega)$ and 
$\eF \in J_{c}(M) = GL_{n}(R)$, the vector field $\eF \eta$ 
belongs to $T[M](\Omega)$. Therefore the relation (\ref{legh}) 
can be written as: 
\begin{equation} 
\frac{\partial^{2} W}{\partial \eF_{ij} \partial 
\eF_{kl}} (F)\ \int_{\Omega} \eta_{i,j} 
\eta_{k,l} \mbox{ d}x \ = \ 0 
\label{leg}
\end{equation}
for any $\eta \in C^{\infty}_{0}$.  An argument from Ball 
\cite{1'}, proof of Theorem 3.4, allow us to consider piecewise 
affine vector fields $\eta$. It can be shown that  (\ref{leg}) 
implies the Legendre-Hadamard equality: 
\begin{equation}
\frac{\partial^{2} W}{\partial \eF_{ij} \partial 
\eF_{kl}} (F)\ a_{i} a_{k} b_{j} b_{l} \ = \ 0  \ \ \ , 
\label{trueleg}
\end{equation}
for any vectors $a,b \in {\mathbb R}^{n}$ . 

Therefore any $W \in nl^{0}(GL_{n})$ is a classical null lagrangian. 

This means that $W$ can be extended over  $R^{n \times n}$ such that 
for any $\eF \in R^{n \times n}$ and for any $\eta \in  C^{\infty}_{0}(\Omega, R^{n})$ 
we have the inequality: 
$$ \int_{\Omega} W(\eF + \nabla \eta) \mbox{ d} x \ = \ \int{\Omega} W( \eF) \mbox{ 
d}x$$
 
The class of classical null lagrangians is known 
(see Ball, Currie \& Olver \cite{21} or Olver \cite{olv}); any 
homogeneous null lagrangian $W$ is a linear combination of minors of $\eF$. From the 
definition of $NL^{0}(GL_{n})$ we see that any $GL_{n}$ null lagrangian is a classical 
null lagrangian. The particular structure of null lagrangians leads to the
introduction of a variational bi-complex in the language of jets.

For the choice $M = SL_{n}$ the relation \eqref{legh} becomes
\begin{equation}
D(DW(\eF) a \otimes b) a \otimes b \ = 0 
\label{hlsln}
\end{equation}
for any orthogonal $a, b \in R^{n}$ (that is $a \cdot b = 0$). 

\begin{prop}
If $W \in nl^{0}(SL_{n})$  then for any $\eF \in SL_{n}$, 
for any $a, b \in R^{n}$, $a \cdot b = 0$ and for 
any $t \in R$ we have $\eF(I_{n} + t a \otimes b) \in SL_{n}$ and the function 
$$t \mapsto f(\eF, a\otimes b, t) =  W(\eF(I_{n} + t a \otimes b) ) $$
is linear. 
\label{pprep}
\end{prop}

\paragraph{Proof.}
Let us denote, for $\eH \in j_{c}(M)$, 
 by $exp_{t} \eH$ the solution of the problem $\dot{\eF_{t}} = \eH
\eF_{t}$, $\eF_{0} = I_{n}$. It is straightforward that if $\eH \eH = 0_{n}$ then 
$exp_{t} \eH \ = \ I_{n} + t \eH$. For any $a, b \in R^{n}$, $a \cdot b = 0$, we 
take $\eH = a \otimes b$ and obtain the first part of the proposition. 

In order to resume the proof, because of theorem \ref{tfundam}, it is 
sufficient to prove  that the second derivative of $ f(\eF, a\otimes b, t) $ 
with respect to $t$  equals $0$ for $t=0$. We use  again that 
if $\eH = a \otimes b$ and $a \cdot b = 0$ then $\eH \eH = 0_{n}$, and we are led 
to the equality: 
$$\frac{d^{2}}{dt^{2}} f(\eF, a\otimes b, t)_{t=0} = D(DW(\eF) a \otimes b) a \otimes 
b$$
This resumes the proof, because of the hypothesis \eqref{hlsln}.
\quad $\blacksquare$

In the case $M = SL_{2}$ we proved the following theorem.

\begin{thm}
If $W \in nl^{0}(SL_{2})$ then $W(\eF) = a_{ij} \eF_{ij} + b$. 
\label{tsl2}
\end{thm}

\paragraph{Proof.}
We shall use local coordinates of $SL_{2}$ and apply proposition \ref{pprep}. 
It is sufficient to consider the coordinates: 
$$ \eF = \left( \begin{array}{cc}
 X & Y \\
Z & \frac{1+YZ}{X} 
\end{array} \right) \ \ , \eF = \left( \begin{array}{cc}
 \frac{1+Y'Z'}{X'} & Y' \\
Z' & X' 
\end{array} \right) $$
Take $a =(a_{1}, a_{2})$ and $b = (-a_{2}, a_{1})$. Then the function 
$f(\eF, a\otimes b, t)$, expressed in the coordinates $(X, Y, Z)$ or $(X', Y', Z')$, 
is linear in $t$, as shown in proposition \ref{pprep}. We derive twice with respect to
time the function $f(\eF, a\otimes b, t)$ at $t=0$ and we equal the result to $0$.
After some elementary computation we obtain the following minimal system of 
equations for 
$g(X,Y, Z) = f(\eF(X,Y,Z))$: 
\begin{equation}
g_{,xx} X^{2} = 2 g_{,yz} (1+ YZ)
\label{5}
\end{equation}
\begin{equation}
g_{,zz} X = - g_{,yz} Y
\label{6}
\end{equation}
\begin{equation}
g_{,xy} X = - g_{,yz} Z
\label{7}
\end{equation}
\begin{equation}
g_{,yy} = 0 
\label{8}
\end{equation}
\begin{equation}
g_{,zz} = 0
\label{9}
\end{equation}
From \eqref{8}, \eqref{9} we find that: 
\begin{equation}
g(X,Y,Z) = A(X) YZ + B(X) Y + C(X) Z + D(X)
\label{10}
\end{equation}
We put the expression of $g$ in \eqref{6} and we obtain the equation: 
$$X C'(X) + XY A'(X) = - A(X) Y$$
From here we derive that $C(X) = c$ and $A(X) = k/X$. We introduce in \eqref{10} 
what we have found and use this in \eqref{7}. We find that $B(X)=b$. Finally, we 
update the form of $g$ and use it in \eqref{5}. It follows that $D''(X) = 2k/X^{3}$ 
therefore $D(X) = (k/X) + eX + f$. We collect all the information gained and 
we arrive at the following expression of $g$:
$$g(X,Y,Z) = k \frac{1+ YZ}{X} + bY + cZ + eX + f$$
which proves the theorem. 
\quad $\blacksquare$ 

A straightforward consequence of the previous theorem is the following one: 

\begin{thm}
Let $W \in C^{2}(SL_{2},R)$ and $\Omega \subset R^{2}$, bounded, open, with 
smooth boundary. If for any volume preserving diffeomorphism $\phi: \Omega 
\rightarrow \Omega$, with compact support we have 
$$\int_{\Omega} W(\nabla \phi) \mbox{ d} x \ = \ ct. $$
then $W(\eF) = a_{ij} \eF_{ij} + b$. 

Moreover, let $W \in C^{2}(GL_{2},R)$. If for any diffeomorphism $\eu: 
\Omega \rightarrow \Omega$ with compact support and for any  volume 
preserving diffeomorphism $\phi: \Omega \rightarrow \Omega$, with compact 
support we have 
$$\int_{\Omega} W(\nabla (\eu .\phi)) \mbox{ d} x \ = \ \int_{\Omega} W(\nabla \eu) 
\mbox{ d}x $$
then $W(\eF) = a_{ij}(\det \eF) \eF_{ij} + b(\det \eF)$. 
\label{trig}
\end{thm}

These results suggest the following conjecture: 

\begin{conje}
For any $M$ and any $W \in nl^{0}(M)$ there is a
classical null lagrangian $\tilde{W}$ such that $W(\eF) = \tilde{W}(\eF)$ 
for any $\eF \in J_{c}(M)$.
\label{con}
\end{conje}

Even if homogeneous $M$ null lagrangians are classical null lagrangians, 
the set $NL^{1}(M)$ can depend non-trivially on $M$. Indeed, consider 
the group of symplectic matrices $M = Sp_{m}$, where $n=2m$. Then 
$[M](\Omega)$ is simply the group of all symplectomorphisms with compact 
support in $\Omega$. Take the one-form $\alpha(x,y) = y \mbox{ d}x$. Then 
$\alpha$ is a  $Sp_{m}$ differential invariant (the Calabi invariant), which give 
raise to a vectorial null lagrangian $\alpha*$, according to proposition \ref{pdinv}. 
We leave the reader to check that (any component of) $\alpha^{*}$   belongs to 
$NL^{1}(Sp_{m})$, but not to $NL^{1}(GL_{n})$.

\section{Final remarks}

It would be very interesting to construct independently the complex of variational 
integrals $VI(M)$. We think that this can be done in the language of currents, 
using a similar approach as  Ambrosio \& Kirchheim \cite{amki}. 

It is to be mentioned that behind the algebraic construction performed in this paper 
are hidden facts related to continuity of variational integrals defined over 
groups of diffeomorphisms, as  in the case of classical null lagrangians 
is shown in Ball, Currie \& Olver \cite{21}. 
In the context of diffeomorphisms groups we cite the lecture paper Buliga \cite{bu}. 
A paper concerning necessary and sufficient conditions 
for lower semicontinuity of variational integrals on diffeomorphisms groups is 
in preparation. 

{\bf Aknowledgements.}Some results from this paper have been communicated in a talk 
given at S.I.S.S.A. Trieste, in September 1999. The author wishes to thank
 A. Braides and G. Dal Maso for the opportunity to give this talk and for interesting 
discussions during the visit.

\end{document}